\documentclass[12pt]{amsart} \usepackage{amssymb}
\usepackage{graphicx}

\theoremstyle{plain}
\newtheorem{theorem}{Theorem}
\newtheorem{lemma}[theorem]{Lemma}

\newtheorem{cor}[theorem]{Corollary}

\theoremstyle{definition}

\begin{document}
%\selectlanguage{angol}

\author{B\'ela Bollob\' as, G\'abor Kun and Imre Leader}

\email{I.Leader@dpmms.cam.ac.uk}
\email{kungabor@cs.elte.hu}
\email{B.Bollobas@dpmms.cam.ac.uk}

\address{Trinity College \\
Cambridge CB2 1TQ \\
UK}
\address{School of Computing Science \\
Simon Fraser University\\
Burnaby BC\\
Canada V5A 1S6}
\address{Trinity College \\
Cambridge CB2 1TQ \\
UK}

\begin{abstract}
We consider the pursuit and evasion game on finite, connected,
undirected graphs known as cops and robbers. Meyniel conjectured
that for every graph on $n$ vertices $O(n^{\frac{1}{2}})$ cops can
win the game. We prove that this holds up to a $log(n)$ factor for
random graphs $G(n,p)$ if $p$ is not very small, and this is close
to be tight unless the graph is very dense. We analyze the
area-defending strategy (used by Aigner in case of planar graphs)
and show examples where it can not be too efficient.
\end{abstract}

\title{Cops and robbers in a random graph DRAFT}

\thanks{Part of this work was supported by OTKA Grant no.
T043671, NK 67867, K67870, by NKTH (National Office for Research
and Technology, Hungary) and by PIMS (Pacific Institute for Mathematical Sciences).}

\maketitle

\section{Introduction}

We will study the following pursuit and evasion game, usually known as
cops and robbers. There is a finite, connected, undirected graph $G$, and $m$ cops and one robber. At the start, each cop chooses one vertex, and then the robber makes his choice of a vertex. Then they move alternately (first the cops then the robber). In the cops' turn, each cop may move to an adjacent vertex, or remain where he is, and similarly for the robber. The cops win the game if one of the cops catches the robber, i.e. lands on the same vertex.
We denote by $c(G)$ the `cop-number' of $G$, meaning the minimal $m$ such that $m$ cops have a winning strategy in $G$, and by $c(n)$ the maximum of $c(G)$ over all graphs with $n$ vertices.

This game has been studied by several authors. Maamoun and Meyniel determined the cop-number for grids \cite{MM}. Aigner and Fromme \cite{AignerFromme}
proved that in the case of planar graphs three cops can catch the robber.
%(and three is actually needed in the edge graph of the icosahedron).
Andreae showed that, for graphs without a complete $K_k$ minor, $k-1
\choose 2$ cops suffice \cite{Andreae}. Quilliot \cite{Q2} found the
upper bound $C(G) \leq 2k+3$ for graphs with orientable genus at
most $k$. Frankl gave lower bounds on $c(G)$ in the case of large
girth graphs \cite{Frankl1}. The graphs with $c(G)=1$ were
characterised by Nowakowski and Winkler \cite{NowakowskiWinkler} and
Quilliot \cite{Q1} independently (such characterizations have
complexity theoretical motivations, see Goldstein and Reingold
\cite{GoldsteinReingold}). We mention in passing that a similar game
was considered by Parsons \cite{Parsons1}, \cite{Parsons2}, in a
continuous setting, but that version is rather different in that the
cops there do not have any information about the robber's moves.

Clearly the most substantial question is to determine the order of magnitude of $c(n)$. Meyniel
conjectured that $c(n)=O(\sqrt{n})$. To see why $\sqrt{n}$ cops can be needed,
note that if a graph $G$ has no cycle of length shorter than five and every vertex has degree at least $\delta$ then $c(G) \geq \delta$: if it is the robber's turn to move then he has a choice
not adjacent to any cop vertex, since each cop has at most one common neighbour with the robber.
In particular we will get $\sqrt{n}$ order of magnitude for the incidence
graph of a finite geometry, that is, the bipartite graph with vertex set consisting of the points and lines of the geometry, and with two vertices representing a point and a line being adjacent if the point is on the line. If the geometry has $q^2+q+1$ points then the bipartite graph will have $2q^2+2q+2$. And at least $q+1$ cops will be needed: at every step the robber will have $q+1$ neighbours and a cop vertex will be adjacent to at most one of these, since the graph contains no triangle and cycle of length $4$ (so the robber will always have an escape choice if the number of cops is at most $q$).

In section 2 we generalize the robber's strategy for large girth graphs. We give a new strategy for the cops.

Our main aim in this paper is to prove that the conjecture essentially holds for sparse random graphs: the cop-number has order of magnitude $\Omega(n^{1/2+o(1)})$ in this case. In fact, our upper bound holds also for denser random graphs, whereas the lower bound does depend on the density. This is the content of Section 3.

The best upper bound known on $c(n)$ is $(1+o(1))\frac{n loglog n}{log n}$, see \cite{Frankl1}. This comes from the simple facts that the neighbourhood of a vertex and also the shortest path between two points can be defended by a single cop. In Section 4, we analyze the question of how efficient such an area-defending strategy can be (with each single cop defending an area independently). It turns out that the area defended by a single cop will be a retract the image -- here by `homomorphism' we mean a mapping of the vertices that
sends each edge to an edge or a single vertex).

Our aim in Section 4 is to prove that such strategy can not be too effective:
we construct a graph $G$ whose largest retract (apart from $G$ itself) has size only a log-power.

Finally, in Section 5 we pose a few open questions.

\section{Srategies}

We will often use the following consequence of Chernoff's theorem about the binomial distribution.

\begin{lemma}~\label{binom}
Let $0 \leq p \leq 1$ and $k, n$ integers, assume $k \leq pn$. Then the inequality

$\sum_{i=0}^k {n \choose i} p^i (1-p)^{n-i} \leq e^{-\frac{(k-pn)^2}{2pn}}$ holds.
\end{lemma}

\subsection{The robber's strategy}

First we give a lower bound on the cop number. The "baby version" of this strategy for the robber was used in large girth graphs by Frankl \cite{Frankl1}. We will assign a weight to every position: this will be some weighted sum of the number of non-backtracking walks of different length from the robber to the cop. The robber will always choose the next step to minimize this function. The robber will move on an induced subgraph $R$ of $G$ with
minimal degree $\delta(R)$ large enough.
Set $M_i(s)=M_i(G,R,s)=max_{x \in V(R), S \subseteq V(G), |S|=s}$ number of non-backtracking walks of length $2i \text{ from } S \text{ to } x.$

\begin{theorem}~\label{condition}
Let $G$ be a connected, graph on $n$ vertices, $R$ an induced subgraph with minimal degree $\delta(R)=d \geq 3$. Then for every positive integer $r$ the inequality $(d-1)^{-r}M_r(2c(G)) > \frac{1}{r+1}$ holds for the cop-number $c(G)$.
\end{theorem}

\begin{proof}
It will be convenient to modify the rules so that the cops are not allowed
to stay at a vertex but all have to move. Note that this is not an important modification, since we need at most twice as many cops to win this game as the original one. Indeed, the cops may go in pairs, with one following the original strategy and the other always going to a neighbour vertex, unless the first cop has to stay according to the original strategy, in which case they swap vertices and swap roles. Hence, in order to prove that $c$ cops cannot catch the robber, it is enough to prove that $2c$ cops always forced to move cannot catch the robber.

Let $N_i$ (depending on the cops' and the robber's position both) denote the number of non-backtracking walks of length $2i$ from the robber to a cop such that the first edge of this walk is not the one the robber used last time. Clearly $N_i \leq M_i(2c)$. We will show that the robber has a strategy against $2c$ cops (forced always to move) if $(d-1)^{-r}M_r(2c) \leq \frac{1}{r+1}$. We will show that the robber has a strategy to keep the following function less than one:

$\begin{displaystyle}
W=\sum_{i=0}^r (d-1)^{i-\frac{i}{r log(d-1)}} N_i.
\end{displaystyle}$

Note that if the robber manages this then he will win, since $N_0 \geq 1$ when the robber is caught.
The robber will also always move and his walk will be non-backtracking. His strategy will be always to minimize $W$.
So assume $W<1$ and that it is the robber's turn. Now according
to the robber's choice a fraction $\geq \frac{d-2}{d-1}$ of the walks in
$W$ is removed from the sum. The cops now make their step:
in the worst case all get closer to the robber, or, to put it
another way, they make the last step of all the walks not
neglected by the robber. Now, a walk of length $2i$ from the robber to a cop corresponds to a walk of length $2(i+1)$ in the previous position.
But we have to be careful with possible backtracking walks in the old position giving non-backtracking walks in the new position. This only can happen in the way that a cop moves from the vertex $x$ to $y$ and so a new type non-backtracking walk $y$ starting with the edges $(yx),(xz)$ for some vertex $z$ contributes to the new sum. But every such walk corresponds to (an even shorter) subwalk from $y$ we counted in the sum but it does not appear, since the cop went in the other way.
Similar thing can not happen on the robber's side: Assume that ther robber moved from $x$ to $y$. Now we may have a walk starting with $(yx)$ which comes from a backtracking walk. But this walk will not contribute to the sum since its first edge is just the one used by the robber last time.

The contribution of
these inherited walks to the new weight $W_{new}$ is at most
$(d-1)^{-\frac{1}{r log(d-1)}}=e^{\frac{-1}{r}}$ times smaller than
it was to $W$: the walk will be shorter by two, but
its weight is $(d-1)^{1-\frac{1}{r log(d-1)}}$ times bigger.
Altogether these give at most $We^{\frac{-1}{r}}$. And we have the
last summand as well: by our assumption this is at most $\frac{1}{r+1} < 1-e^{\frac{-1}{r}}$.
This yields $W_{new} \leq (1-e^{\frac{-1}{r}})W+e^{\frac{-1}{r}}$. Hence if
$W<1$ then $W_{new}<1$.

Finally we have to find an appropriate initial position with $W<1$.
(Here we modify the rules and assume that the robber will also choose the initial position for the cops and makes the first step: this makes no difference.)

We will choose the initial position of the cops and the robber randomly (according to the uniform distribution) and we prove that the expected value of $W$ is less than one. Let $E_i$ denote the expected number of non-backtracking walks from a vertex $x \in V(R)$ to a set $S \subseteq V(G)$ of size $2c$. Clearly $E_i \leq M_i(2c)$. Note that $E_{i+1} \geq (d-1)^2 E_i$: the number of non-backtracking walks of length $2(i+1)$ from $R$ is at least the number of non-backtracking walks of length $2i$ from $R$ with the first two steps in $R$ and this is at least $(d-1)^2$ times the number of walks of length $2i$ from $R$, since $\delta(R)=d$. Now the expected value
of $W$ is

$\sum_{i=0}^{r} (d-1)^{-i} E_i < (r+1) (d-1)^{-r} E_r \leq (r+1) (d-1)^{-r} M_r(2c) < 1$.

This completes the proof of the theorem.
\end{proof}

\subsection{The cops' strategy}

Our aim is to place the cops on some `spread-out' set of vertices,
and the hope is then that, wherever the robber may be,
our cops are dense enough near a ball around the robber that they can move in and seal him off. This `dense enough' will be accomplished by a Hall type of argument.

For $x \in V(G)$ and a positive integer $r$ we denote by $B(x,r)$ the ball of radius $r$ around
$x$: the set of vertices at distance at most $r$ from $x$. For a set of vertices $S$, let $N(S,r)$ denote the $r$-neighbourhood of $S$, that is, the set of vertices at distance at most $r$ from $S$.

\begin{lemma}
Let $G$ be a connected graph on $n$ vertices, and let $I \subseteq V(G)$. If for every $x \in V(G)$ there exists $r=r_x$ such that for every $S \subseteq B(x,r)$ the inequality $|I \cap N(S,r+1)| \geq |S|$ holds, then $|I| \geq c(G)$.
\end{lemma}

\begin{proof}
We give a winning strategy with $|I|$ cops. The cops' initial position is the set $I$.
Let $x$ denote the robber's vertex and $r$ the corresponding radius (as given in the
statement of the Lemma). By the K\H onig-Hall theorem
we can assign to every $y \in B(x,r)$ a cop in $B(y,r+1)$ such that we assign to
every vertex $y$ a different cop. So let each cop initially, in the first $(r+1)$ steps, go to his vertex and stay there. Since the robber cannot leave $B(x,r)$ in $r$ steps he will be caught.
\end{proof}

\begin{cor}~\label{cop}
Let $G$ be a connected graph on $n$ vertices and $c$ an integer. If there exists an integer $r$ such that for every $x \in V(G), S \subseteq B(x,r)$ the inequality $10|S|log(n) \leq \frac{c}{n}|N(S,r+1)|$ holds, then $c \geq c(G)$.
\end{cor}

\begin{proof}
The case $n \leq 3$ is trivial, so we assume $n \geq 4$.
We will choose a random initial position $I$ for the cops and show that the conditions of the previous lemma will hold with positive probability. For every $x \in V(G)$ the probability that $x \in I$ will be $\frac{c}{2n}$ and these events will be independent. The probability that $|I|>c$ is less than $1/2$ by the Markov inequality. Consider the vertex $x \in V(G)$ and the subset $S \subseteq B(x,r)$. The probability that $|N(S,r+1) \cap I| < |S|$ is at most $e^{-\frac{(|S|-c|N(S,r+1)|/2n)^2}{c|N(S,r+1)|/n}}$ by Lemma~\ref{binom}. We know that $c|N(S,r+1)|/2n-|S| \geq \frac{2}{5} c|N(S,r+1)|/n \geq 4|S|log(n)$ since $10|S|log(n) \leq \frac{c}{n}|N(S,r+1)|$.
Altogether the probability of the event that $|N(S,r+1) \cap I| < |S|$ is less than $e^{-\frac{8}{5}|S|log(n)}=n^{-\frac{8}{5}|S|}$. The probability that such an $S$ exists is at most $\sum_{s=1}^{\infty} {n \choose s} n^{-\frac{8}{5}s} \leq \sum_{s=1}^{\infty} n^{-\frac{3}{5}s}<1/2$ since $n \geq 4$. We know that $|I| \leq c$ holds with probability at least $1/2$, hence there exists an appropriate $I$ of size $c$.
\end{proof}

\section{Random graphs}

Let $G=G(n,p)$ denote the Erd\H os-R\'enyi random graph. We will estimate $M_r$ in this graph model in order to get a lower bound on the cop-number. First we need to estimate the size of small balls.

\begin{lemma}
With high probability the following holds for the random graph $G(n,p)$ if
$pn>1$: for every integer $k$ and $x \in V(G)$ the ball $B(x,k)$ has size at most $20 log(n) (1+pn)^k$.
\end{lemma}

\begin{proof}
Given a set $S \subseteq V(G)$ the probability that the size of $N(S,1)=N(S)$ differs by at least $\sqrt{3pn|S|log(n)}$ from its expected value is at most $n^{-\frac{3|S|}{2}}$ by Lemma~\ref{binom}. So the probability that such a set exists is at most $\sum_{s=1}^{\infty} {n \choose s} n^{-\frac{3s}{2}}<\sum_{s=1}^{\infty} n^{-\frac{s}{2}}=o_n(1)$. Hence whp this holds for every $S$: we will assume this in what follows. For every set $S$ the expected size of $N(S,1)$ is at most $(1+pn)|S|$. Hence we have $|N(S,1)| \leq (pn+1)|S|+\sqrt{3|S|pnlog(n)}$. So $|B(x,1)| \leq pn log(n)$ for every $x \in V(G)$ if $n$ is large enough. We will prove that $|B(x,k)| \leq log(n)(pn+1)^k e^{\sqrt{3} \sum_{i=2}^k 2^{-i/2}}$ by induction on $k$. This holds for $k=1$. Assuming the statement for $k$ we have

$|B(x,k+1)| \leq |B(x,k)|(pn+\sqrt{\frac{3pnlog(n)}{|B(x,k)|}}) \leq \\
\Big( log(n)(pn+1)^k e^{\sqrt{3} \sum_{i=2}^k 2^{-i/2}} \Big) \Big(pn+1+\sqrt{\frac{3pnlog(n)}{log(n)(pn+1)^k e^{\sqrt{3} \sum_{i=2}^k 2^{-i/2}} }} \Big) \leq \\
\Big(log(n)(pn+1)^k e^{\sqrt{3} \sum_{i=2}^k 2^{-i/2}} \Big) (pn+1) (1+\sqrt{\frac{3}{(pn+1)^{k+1}}}) \leq \\ \Big( log(n)(pn+1)^k e^{\sqrt{3} \sum_{i=2}^k 2^{-k/2}} \Big) (pn+1) (1+e^{\sqrt{3}2^{-\frac{k+1}{2}}}) \leq \\ log(n)(pn+1)^{k+1} e^{\sqrt{3} \sum_{i=2}^{k+1} 2^{-i/2}}$.

Since $e^{\sqrt{3} \sum_{i=2}^{\infty} 2^{-i/2}}= e^{\sqrt{3}\frac{2+\sqrt{2}}{2}}<20$ the lemma follows.
\end{proof}

\begin{lemma}
With high probability the following holds for the random graph $G(n,p)$ if
$pn>1$: for every $0<\varepsilon<\frac{1}{3}$, integers $k$ and

\noindent
$r<\frac{((1/2-\varepsilon)log(n)-loglog(n)-log(40))}{log(pn+1)}-1$
and for every pair of vertices $x, y \in B(x,r)$ the number of non-backtracking paths of length $k$ from $x$ to $y$ in $B(x,r)$ is
at most $(\frac{7}{\varepsilon})^k$.
\end{lemma}

\begin{proof}
We have seen that whp for every $x \in V(G)$ the ball $B(x,r)$ has size
at most $\frac{n^{\frac{1}{2}-\varepsilon}}{2(pn+1)}$. Given a set $S$ the expected number of the edges with at least one endpoint in $S$ is at most $pn|S|$. Lemma~\ref{binom} shows that the probability that the number of such edges differs by at least $\sqrt{3pn|S|log(n)}$ from its expected value is at most $n^{-\frac{3|S|}{2}}$. So whp this does not
hold for any $S \subseteq V(G)$ since $\sum_{s=1}^{\infty} {n \choose s} n^{-\frac{3s}{2}}= O(n^{\frac{3}{2}})$. In particular, we may assume for every ball $B(x,r)$ that it has at most $n^{1/2-\varepsilon}$ edges: the expected value is at most $pn|B(x,r)| \leq \frac{1}{2} n^{1/2-\varepsilon}$ and $\sqrt{3pn|B(x,r)|log(n)} \leq \frac{1}{2}n^{1/2-\varepsilon}$ if $n$ is
large enough.

Now we prove that whp for every vertex $x \in V(G)$ the ball $B(x,r)$ consists of a tree plus at most $\frac{3}{\varepsilon}$ edges. We may think about $B(x,r)$ as a set defined in a process
of $r$ steps: We start with $\{x\}$, then we add the neighbours of $x$, the new neighbours of this set etc. For every new vertex and new edge from this point the probability that the other endpoint of this edge is an old point is at most $n^{-1/2-\varepsilon}$.
Hence the probability to have at least $3/\varepsilon$ such edges is at most ${n^{1/2-\varepsilon} \choose 3/\varepsilon} (n^{-1/2-\varepsilon})^{3/\varepsilon}<n^{-\frac{3}{2}}$. Whp the number of such edges is at most $3/\varepsilon$ in every ball of radius $r$. After the removal of these edges from the ball we get a tree.
Every non-backtracking path in this ball is completely described by its endpoints and the used edges (with direction) not in the tree. Hence the number of such paths of length $k$ is at most $(1+6/\varepsilon)^k < (\frac{7}{\varepsilon})^k$.
\end{proof}

\begin{lemma}
With high probability the following holds for the random graph $G(n,p)$ if
$pn>1$: for every $0<\varepsilon<\frac{1}{3}$, integer

\noindent
$r<\frac{(1/2-\varepsilon)log(n)-loglog(n)-log(40)}{log(pn+1)}-1$
and for every pair of vertices $x, y \in V(G)$ the number of non-backtracking paths of length $\leq 2r$ from $x$ to $y$ is at most $(\frac{7}{\varepsilon})^{3r}$.
\end{lemma}

\begin{proof}
We know by the previous lemma that whp for every $z \in B(x,r)$ the number of non-backtracking paths of length $k$ from $x$ to $z$ in $B(x,r)$ is at most $(\frac{7}{\varepsilon})^k$, and the same holds for the ball $B(y,r)$ and a vertex $z \in B(y,r)$. First we estimate the number of those paths $x_0=x, \dots ,x_l=y$, where for the last vertex $x_k \notin B(y,r)$ of the path
$x_k, x_{k+1} \notin B(x,k-1)$ holds. We call such paths \emph{special}.

Set $S=V(G) \setminus (B(x,k-1) \cup B(y,r)).$ Consider a permutation $\nu$ of $S$ and the graph $G_{\nu}$ with the vertex set of $G$ and edge set $E(G_{\nu})=E(G) \cup \{(a,b): a \in B(y,r) \setminus B(x,k-1), b \in S, (a, \nu(b)) \in E(G) \} \setminus \{(a,b): a \in B(y,r) \setminus B(x,k-1), b \in S, (a, b) \in E(G) \}$. Note that the graphs $G(n,p)=G_{\nu}$ are equiprobable.
The edges from $B(y,r)\setminus S(x,k-1)$ to $S(x,k)$ are the ones where these graphs $G_{\nu}$ may disagree. E.g. $(x_k,x_{k+1})$ is such an edge for a special path.

Now we examine the set
of edges from the sphere $S(x,k)$ to $B(y,r) \setminus B(x,k-1)$. We know that whp $B(y,r)$ and $S(x,k)$ have at most $n^{\frac{1}{2}-\varepsilon}$ vertices and edges for every $x,y \in V(G)$. And the distribution of the induced subgraph on the vertices $V(G) \setminus B(y,r-1)$ has the same distribution as $G(n-|B(y,r-1)|,p)$.
The probability for an edge with one endpoint in $S(x,k)$ that the other endpoint is in $B(y,r)$ will be at most $\frac{n^{1/2-\varepsilon}}{n-3n^{1/2-\varepsilon}}=\frac{1}{n^{1/2+\varepsilon}-3}$. The probability to have at least $\frac{2}{\varepsilon}$ edges joining $S(x,k)$ and $B(y,r)$ is at most ${n^{1/2-\varepsilon} \choose \frac{2}{\varepsilon}} (n^{1/2+\varepsilon}-3)^{\frac{2}{\varepsilon}}=o(n^{-2}).$
So whp this does not hold for any $x,y \in V(G)$. The number of special paths (where $x_k, x_{k+1} \notin B(y,r)$) is at most $\sum_{k=1}^{2r} \frac{2}{\varepsilon} (\frac{7}{\varepsilon})^k (\frac{7}{\varepsilon})^{2r-k-1}<\frac{r}{2}(\frac{7}{\varepsilon})^{2r}$.

Now we estimate the total number of paths from $x$ to $y$. Consider a path $x_0=x, \dots ,x_l=y$ such that $l \leq 2r$, and the first point in $B(y,r),x_{k+1}$ is in $B(x, k-1)$. Let $a_0=x, \dots ,a_j=x_{k+1}$ a shortest path from $x$ to $x_{k+1}$ and consider the path $y_0=a_0=x, y_1=a_1, \dots ,y_j=a_j=x_{k+1}, y_{j+1}=x_{k+2}, \dots y_{l-j+k-1}=x_l=y$. Note that this path is special.
To every such special path we did correspond at most $(\frac{7}{\varepsilon})^{k+1}$ paths, since this is an upper bound for the number of paths from $x$ to $x_{k+1}$ with length $(k+1)$. Since $k+1 \leq r$ the number of such paths is at most $r(\frac{7}{\varepsilon})^{2r}\sum_{k=0}^r (\frac{7}{\varepsilon})^{k+1}<2r (\frac{7}{\varepsilon})^{3r}$.

Now consider a path such that $x_{k+1} \notin B(x,k-1)$ but $x_k \in B(x,k-1)$.
To such a path we can correspond again a special one by replacing the subpath $x=x_0, \dots ,x_k$ by a shortest path from $x$ to $x_k$ (with length $(k-1)$).
We correspond to every special path at most $(\frac{7}{\varepsilon})^k$ paths.
By the same argument as above we get that the number of such paths is at most
$\frac{r}{2}(\frac{7}{\varepsilon})^{3r-1}$.

Finally, in the case $y \in B(x,r)$
there are paths from $x$ to $y$ completely inside $B(y,r)$. The number of these paths is at most $\sum_{l=0}^{2r} (\frac{7}{\varepsilon})^l < 2 (\frac{7}{\varepsilon})^{2r}$. Altogether, suming up the number of four types of path we get that whp for every $x,y \in V(G)$ there are at most $r(\frac{7}{\varepsilon})^{3r}$ non-backtracking paths with length at most $2r$ from $x$ to $y$.
\end{proof}

\begin{theorem}
The following lower bound holds for the cop-number of $G(n,p)$ with probability going to $1$ as $pn \rightarrow \infty$: \\

$c(G)>\frac{1}{(pn)^2} n^{\frac{1}{2} \frac{loglog(pn)-9}{loglog(pn)}}$.
\end{theorem}

\begin{proof}
First we will find a nonempty induced subgraph $R$ of $G$ with minimal degree at least $\frac{pn}{4}$.
All but $\frac{n}{5}$ vertices have degree at least $\frac{3}{4}pn$ with high probability if $pn$ is large enough. Let $B$ denote the small set of these exceptional vertices. Consider the maximal set of vertices $R$ with the following properties:

\begin{enumerate}

\item
$B \cap R = \emptyset$,

\item
Every $x \in R$ has more than $\frac{1}{4}pn$ neighbours in $R$.

\end{enumerate}

We show that $|V(G)|-|R| \leq 4|B| \leq \frac{4n}{5}$, hence $R \neq \emptyset$. The set $R$ will not contain the vertices of $B$, those vertices with too many neighbours in $B$, and those with too many neighbours in this set etc. We may think about the definition of the complement of $R$ as a process, where we decide about new and new vertices with too many neighbours in $R^c$ to be in $R^c$. When we decide that a vertex $x$ is in $R^c$ then $x$ have at least $\frac{1}{2}pn$ neighbours decided to be in $R^c$ and at most $\frac{1}{4}pn$ other neighbours. So the edge boundary of the points in $R^c$ decreases by at least $(\frac{pn}{4}+1)$ when adding such a point, and the edge boundary of $B$ was at most $(\frac{3}{4}pn+3)|B|$.

We know $\delta(R) \geq pn/4$. Set
$r=[\frac{(1/2-\varepsilon)log(n)-loglog(n)-log(40)}{log(pn+1)}]-1$.
and $\varepsilon=\frac{4}{loglog(pn/4)}$.
We know by the previous lemma that $M_r(1)<r(\frac{7}{\varepsilon})^{3r}$.
Clearly $M_r(2c) \leq 2c M_r(1)$.

Theorem~\ref{condition} yields that $\frac{1}{r+1}<(pn/4-1)^{-r}M_r(2c(G)) \leq 2c(G)(pn/4-1)^{-r}M_r(1)$. Hence

\noindent
$c(G) > \frac{1}{2(r^2+r)}
(\frac{64(pn/4-1)}{343 loglog(pn/4-1)})^r > \\
\frac{1}{2(r^2+r)} (\frac{16pn-64}{343 loglog(pn/4-1)})^{\frac{(1/2-\varepsilon)log(n)-loglog(n)-log(40)}{log(pn+1)}-2}
>\frac{1}{(pn)^2} n^{\frac{loglog(pn)-8+o_{pn}(1)}{2loglog(pn)}}$.

And this is greater than $\frac{1}{(pn)^2} n^{\frac{loglog(pn)-9}{2loglog(pn)}}$ if $pn$ is large enough.
\end{proof}

%%%%%%%%%%%%%%%%%%%%%%%%%%%%%%%%%%%%%%%%%%%%%%%%%%%%%%%%%%%%%%%%%%%%%%%%%%%%%%%%%%%%%%%%%%%%%%%%%

Now we will prove an upper bound on the cop-number.
We will estimate the vertex expansion in random graphs in order to use Corollary~\ref{cop}.

\begin{theorem}
Let $0<\varepsilon<1$. With high probability the following upper
bound holds for the cop-number of the random graph $G=G(n,p)$ if
$p>2(1+\varepsilon)log(n)/n$:

\noindent
$c(G(n,p))<n^{\frac{1}{2}}log(n)max \{ \frac{1}{\varepsilon}; 160000 \}$.
\end{theorem}

\begin{proof}
Consider a subset $S \subset V(G)$. The expected size of the vertex neighbourhood $|N(S)|$ is $|S|+(n-|S|)(1-(1-p)^{|S|}))$. This is less than $(pn+1)|S|$. On the other hand it is at least $n(1-e^{-p|S|})$. Given a subset $S \subseteq V(G)$ the probability that $|N(S)|$ differs by at most $\sqrt{2log(n)(1+\varepsilon/2)(pn+1)|S|} \leq \frac{4-\varepsilon}{4}(pn+1) \sqrt{|S|}$ from its expected value is at least $1-n^{-\frac{2+\varepsilon}{2}|S|}$ by Lemma~\ref{binom}. Whp this holds for every $S \subseteq V(G)$. We will assume this in what follows for every $S$.

Set $r=[\frac{1000log(n^{\frac{1}{2}})}{log(pn+2)}]$.
In order to get a lower bound on the vertex expansion we will find some upper bound first to use our condition conveniently.
If $|S| \geq (pn+1)^2$ then $|N(S)| \leq (pn+2)|S|$. And for every set $S$
we have the inequality $N(S) \leq (pn+1)|S|+(pn+1)\sqrt{|S|} \leq 2(pn+2)|S|$. Now by induction for every integer $k$ and every $S \subseteq V(G)$ the inequality $N|(S,k)| \leq 4 (pn+2)^k$ holds.
In particular, $|B(x,r)| \leq 4000 n^{\frac{1}{2}}$ for
every $x \in V(G)$.

We will estimate for every $x \in V(G)$ and every subset $S \subseteq B(x,r)$ the size of $N(S,r+1)$. We will succeed by induction estimating $N(S,k)$ for every $k \leq r+1$. We must be careful in the cases $k=1,2$ and $k=r,r+1$, this makes our calculations quite technical looking.

The expected value of $|N(T)|$ is $>n(1-e^{-p|T|})$ for every $T$. If $T$ is not too large, namely $p|T|(log(n)+1)<1$ then this is at least $(1-\frac{1}{log(n)})pn|T|$. This holds for $T=N(S,k)$ if $S \subseteq B(x,r)$ for some $x \in V(G)$ and $k \leq r-2$ (these imply $N(S,k) \leq \sqrt{n}(pn+2)^k$).

First, $|N(S)| \geq \frac{\varepsilon}{5} pn|S|$ if $n$ is large enough and
$p|S|(log(n)+1) \leq 1$. Secondly, if $|N(S)|p(log(n)+1)<1$ then $|N(S,2)|
\geq \frac{\varepsilon}{6} p^2n^2|S|$ if $n$ is large enough.
If $3 \leq k \leq r-2$ then the difference from the expectation can not be significant: $|N(S,k+1)| \geq (1-\frac{1}{log(n)})(pn-1)|N(S,k)|$.
Using this for $k=2, \dots ,(r-2)$ we get

\noindent
$|N(S,r-1)| \geq \frac{\varepsilon}{6}(1-1/log(n))^{r-1} (\frac{pn-1}{pn+2})^{r-1} (pn+2)^{r-1} |S|=\\
(1+o(1))\frac{\varepsilon}{6} (pn+2)^{r-1} |S|$.

Next, $|N(S,r)| \geq n(1-e^{-p|N(S,r-1)|})$. This is at least $n/4$ if $p|N(S,r)|>1/2$ and else at least $(2-2e^{-1/2})pn|N(S,r-1)|>1/2 (pn+2)|N(S,r)|$
(assuming again that $n$ is large enough). Altogether we get that $|N(S,r)| \geq max \{n/4; 1/2(pn+2)|N(S,r-1)| \}$. Applying the same argument again $|N(S,r+1)| \geq max \{n/4; 1/4(pn+2)|N(S,r)| \} \geq max \{ n/4; \frac{\varepsilon}{96}(1+o(1))(pn+2)^{r+1}|S| \}$ follows. This is at least
$max \{ n/4; \frac{\varepsilon}{100}(1+o(1))(pn+2)^{r+1}|S| \}$ if $n$ is large enough. By the choice of $r$ we have $\frac{\varepsilon}{100} (pn+2)^{r+1} \geq 10 \varepsilon \sqrt{n}$.

Corollary~\ref{cop} implies that $c(G)< 10 n log(n) max_{x \in V(G), S \subseteq B(x,r)}
\frac{|S|}{|N(S,r+1)|} \\
\leq 10 n log(n) max \{ \frac{1}{10 \varepsilon n^{\frac{1}{2}}}; \frac{4000 n^{\frac{1}{2}}}{n/4} \} = log(n) n^{\frac{1}{2}} max \{ \frac{1}{\varepsilon}; 160000 \}$.

\end{proof}

%%%%%%%%%%%%%%%%%%%%%%%%%%%%%%%%%%%%%%%%%%%%%%%%%%%%%%%%%%%%%%%%%%%%%%%%%%%%%%%%%%%%%%%%%%%%%%%%%%

\section{Defending an area}

In this section we analyze the area-defending strategy. By this strategy we
mean that every single cop defends an area by himself, where `defends an area'
means `moves around in that area in such a way that, if the robber ever
enters the area, he is instantly caught by the cop'. A moment's thought
shows that the
area-defending-strategy of a cop is a retraction $r: G \rightarrow
G$, that is, a homomorphism of the reflexive graph (i.e. the image
of an edge is either an edge or a single vertex) which fixes its
image: $r \circ r=r$. When the robber is at the vertex $x \in G$
then the cop goes to $r(x)$.

We prove that this strategy can not
be too effective: in some graphs the largest area that can be defended by one
cop (equivalently, the largest image of a non-identity retract) is at most
a power of log.

\begin{theorem}
For every positive integer $n$ there is a graph on $n$ vertices
with largest proper retract of size $O(log(n)^8)$.
\end{theorem}

\begin{proof}
First we choose three positive integers $d,s,t$: the choice will
depend on $n$. These will satisfy the conditions $l>s+8$ and
$2ds>(2d-1)(l+2)$. We specify the other conditions on $d,s$ and
$t$ at the end of the proof: the precise values are important only to ensure that we obtain a graph on exactly $n$ vertices.

Consider the $d$-dimensional hypercube $Q$. We subdivide every
edge of $Q$ by adding either $s, (s+1), (s+2)$ or $l, (l+1)$ or
$(l+2)$ vertices. We call the edges of $Q$ divided by $s, (s+1)$
or $(s+2)$ vertices `short' and the other edges `long'. We decide for
every edge independently and randomly if the edge is short or long
(each with probability one half). This random choice gives us
many graphs, because every long edge may have $(l+1), (l+2)$ or
$(l+3)$ edges and every short one $(s+1), (s+2)$ or $(s+3)$.
We will prove that with high probability none of these graphs will
have a large proper retract. And finally we will choose the lengths
of the edges to have a graph on exactly $n$ vertices.

We denote the resulting graph by $G$. We call a subgraph $S$ of
$G$ {\it reduced} if for every $x \in S$ either $x \in V(Q)$ or
the complete path corresponding to this edge of $Q$ containing $x$
is in $S$, and if $S$ contains the endpoints of an edge of $Q$
then the corresponding path is in $S$. So a reeuced graph $S$ is
determined by $S \cap Q$. We denote the following
subgraph of $Q$ by $S'$: $V(S')=V(S) \cap V(Q)$ and $x,y \in V(Q)$
are adjacent in $S'$ if they are adjacent in $Q$ and $S$ contains
the path connecting them. Every retract $R$ gives rise to a reduced retract
of size $\geq [\frac{|R|}{d(l+2)}]$: $R$ is connected
and so every "bad" vertex of $R$ is on an edge of $Q$, and can be
mapped to the endpoint which is in $R$.

We will show that with positive probability all reduced retracts of
$G$ are of size $O(ld^5)$. First we show that the number of the
corresponding sets $R'$ is small. Consider the vertices of
$Q$ as $0-1$ vectors of length $d$, with the $i$th coordinate of $x \in
V(Q)$ denoted by $x_i$. We say that the subset $S
\subseteq V(Q)$ is a {\it union of quarters} if for all $x \in
V(S)$ there are two coordinates $1 \leq i < j \leq d$ such that if
for a vertex $y \in V(Q)$ we have $x_i=y_i$ and $x_j=y_j$ then $y
\in V(S)$. The number of such subsets of $Q$ is clearly at most
$2^{2d^2-2d}$.

\medskip

{\bf Claim:}
Let $r:G \rightarrow G$ be a retraction and $R=r(G)$ a reduced retract of $G$. Then $V(Q) \setminus V(R')$ is a union of quarters.

\medskip

First assume that there is a $z \in V(Q) \cap V(R')$ adjacent to $x$. We have $dist_G(x,z) \geq
dist_G(r(x),z)$. Let $v$ denote the closest vertex to $r(x)$ in $Q$. We know that either $v=z$ or else $v$ is adjacent to $z$, in which $j$ denote the coordinate where they differ. (If $v=y$ then any coordinate will suffice as $j$.) The vertices $x$ and $z$ differ in coordinate $i$, where $j \neq i$. We will show that if $y \in Q$ agrees with $x$ in coordinates $i,j$ then $y \notin R'$. It will
follow that $dist_G(x,y)<dist_G(r(x),y)$. We know that $dist_G(r(y),r(x)) \geq dist_G(y,r(x)) \geq dist_G(y,v)-\frac{l+2}{2} \geq dist_Q(y,v)(s+1)+\frac{l+2}{2} >
dist_Q(y,x)(l+2) \geq dist_G(y,x)$.

If $x$ is not adjacent (in $Q$) to any vertex in $R'$ then consider a shortest path (in $Q$) from $x$ to $R'$. Let $z$ denote its endpoint and $x' \in Q$ the vertex of the path adjacent to $z$ (in $Q$). Now we know that $x'$ is contained by an appropriate "quarter": if $x'$ and
$z$ differ in coordinate $i$
and $r(x')$ is at distance at most one from the edge corresponding to coordinate $j$ this quarter will correspond to coordinates $i,j$. The vertex $x$ has to agree with $x'$ in these coordinates, as otherwise there would be a shorter path from $x$ to $R'$ in $Q$ and so in $G$.
So $V(Q) \setminus V(R')$ is a union of quarters.

\medskip

Now we will show that with high probability $G$ has no large proper reduced retract. We will show that for every large induced subgraph $R'$ whose complement is a union of quarters the probability that the (unique) reduced subgraph $R$ which $R'$ corresponds to is a retract of
$G$ is small.

%In what follows $R$ will denote a potential reduced retract, this $R \cap Q$ is the intersection of quarters in the above sense, and for every $x \in R$ either $x \in V(Q)$ or the whole path corresponding to this edge of $Q$
%is in $R$, moreover if $R$ contains two adjacent vertices in $Q$ then it also will contain the path between these.

We will use the following two observations.

\begin{enumerate}

\item Let $x_1, x_2, x_3, x_4$ be a 4-cycle in $Q$. If a reduced retract of $G$ contains
$x_2, x_3, x_4$ but not $x_1$ then $ dist_G(x_2,x_3)+dist_G(x_3,x_4) \leq dist_G(x_2,x_1)+dist_G(x_1,x_4)$.
The probability that this event is possible, i.e. there are not more long edges on the left hand side than on the right hand side, is $\frac{11}{16}$.

\item Let $x_1, x_2 ,x_3, x_4, y_1, y_2, y_3, y_4 \in Q$ be a ($3$ dimensional) subcube of $Q$, where both $x_1, x_2, x_3, x_4$ and $y_1, y_2, y_3, y_4$ are
4-cycles. If a reduced retract contains $x_1, x_2, x_3, x_4$ but not $y_1, y_2, y_3, y_4$ then $dist_G(y_1,y_2)+dist_G(y_2,y_3)+dist_G(y_3,y_4)+ dist_G(y_4,y_1) \geq dist_G(x_1,x_2)+dist_G(x_2,x_3)+dist_G(x_3,x_4)+dist_G(x_4,x_1)$. The probability that this event is possible is $\frac{163}{256}$.
\end{enumerate}

In the first case the distance of $x_2$ and $x_4$ in $R$ is the length of the path via $x_3$ connecting them, and this cannot be longer then the path via $x_1$ if $R$ is a retract. It is easy to check that the probability of this event is $\frac{11}{16}$. In the second case the cycle of length four via the points not in the retract has to be shorter than the one in the retract: this needs some case analysis. We only will use the fact that the probability of at least one event occurring strictly less than one.

We will show that for every adjacent pair $(a,b)$ (in $Q$), where  $a \in R', b \in Q \setminus R'$  either there are four vertices $y_1 \in Q \setminus R , y_2, y_3, y_4 \in Q \cap R$ at distance at most $3$ from $x$ forming the first configuration or there are eight vertices at distance at most $3$ forming the second configuration.

We may assume without the loss of generality that $a$ and $b$ differ in coordinate $1$. Consider $b_1 \in R'$ adjacent to $b$ (say these differ in the second coordinate) and $a_1 \in Q$ adjacent to both $b_1$ and $a$. If $a_1 \in R'$ then these give an appropriate cycle: $x_1=b, x_2=a, x_3=b_1, x_4=a_1$.
Let us assume $a_1 \notin R'$. Now we pick another vertex $b_2 \in R'$ adjacent to $b$ or $b_1$: say it is adjacent to $b_1$ and the third coordinates differ. Let $b_3$ denote the common neighbour of $b$ and $b_2$. We may assume that $b_2 \notin R'$, as otherwise $b, b_1, b_2, b_3$ would form a cycle we are looking for. Now assume that the vertex $a_3$ the common neighbour of $b_3$ and $a$ is not in $R'$, as otherwise $a, b, b_3, a_3$ would form an appropriate cycle. Similarly we may assume that the vertex $a_2$ which is
the common neighbour of $a_1$ and $b_2$ is in $R$. These eight points form the second configuration:
$x_1=a_1, x_2=a_2, x_3=a_3, x_4=a \notin R'$ and $y_1=b_1, y_2=b_2, y_3=b+3, y_4=b \in R'$.

%(Note that
%none of these edges corresponds to coordinate $i$, else $x,y,y_2$ is a required configuration.)
%We may assume that the fourth vertex of the $4$-cycle containing $y, y_1, y_2$ is also in $R$ else
%we could find again the first configuration. Now either the other face of this square (in direction
%i$) is not in $R$ or we could find the first configuration in this hypercube.

We know that the number of induced subgraphs of $S \subseteq Q$
such that $V(Q) \setminus V(S)$ is a union of quarters is at most
$2^{2d^2-2d}$. We will show that for every such potential $R'$
large enough there are many such disjoint configurations. We know
by Harper's theorem that the edge boundary of a subset of the
hypercube is at least the size of the subset (if the subset has
size at most half of the hypercube). The complement of $S$ has
size at least $\frac{1}{4}|Q|$ unless $S=Q$, so the boundary has
at least $min \{|S|, |Q|/4 \}$ vertices. Let us find a bad
configuration in $Q$, then a new bad configuration not covered by
the $3$-neighbourhood of this configuration, and so on. The
$3$-neighbourhood of a bad configuration has size at most $8 (
{d-3 \choose 3} + {d-3 \choose 2} + {d-3 \choose 1} + 1 ) <
\frac{4}{3} d^3$. So if $min \{ |Q|/4,|S| \} > 3 d^5 > 4/3 d^3
\frac{2d^2 log(2)}{log(16/5)}$ then the probability that $S=R'$
for a reduced retract $R$ of $Q$ is at most $2^{-2d^2+2d-1}$. So
the probability that $G$ has no retract of size at least
$3(l+2)^2d^6$ is at least $\frac{1}{2}$.

Now we specify our other conditions on $s,l$ and $d$. We need
$|d2^{d-1} \frac{s+l+2}{2}-n|+2(l-s)\sqrt{2d2^{d-1}} < d2^{d-l}$:
with probability $\frac{1}{2}$ the difference between the number of
long and short edges is at most $2(l-s)\sqrt{2d2^{d-1}}$ by the
Chebyshev inequality. So the left hand side will be at most the
number of vertices of $G$ (assuming that short edges have $s+2$
edges and long edges have $l+2$ edges all) minus $n$. The right
hand side is the number of edges in $Q$: at every edge we can add
or remove a vertex to have exactly $n$ vertices. So the conditions
are: $l>s+8, 2ds>2(d-1)(l+2)$ and $|d2^{d-1}
\frac{s+l+2}{2}-n|+2(l-s)\sqrt{2d2^{d-1}} < d2^{d-l}$. Let us
choose $d$ such that $d2^{d-1} 100d \leq n \leq (d+1)2^d
100(d+1)$, and let $s+l+2$ be the closest integer to $\frac{n}{d2^d}$
and $l=s+9$ or $l=s+10$. This choice will satisfy the conditions
if $n$ (and so $d$) is large enough.

\end{proof}

There are many nice topological techniques to show that there is no
homomorphism from one graph to another.
%in particular that a graph has large chromatic number. Retractions of reflexive graphs allow to map an edge to a single vertex, even the whole graph may map to any single vertex. The theorem deals with a mixed problem: an extremal/algebraic (topological) problem.
We would be glad to see some interesting, say topological, example. The
graph constructed in the proof is quite close to a product, so
similar strategies will work in this graph like in product graphs
(or particularly grids).

\section{Open questions}

We start by repeating what must be the main open question:

{\bf Question 1.} What is the order of magnitude of the function $c(n)$?

Recall that Meyniel conjectured $c(n)=O(\sqrt{n})$. The best upper bound on
$c(n)$ is that of Frankl \cite{Frankl1}, namely $(1+o(1))\frac{n loglog n}{log n}$.
Thus even an upper bound of $n^{1-\varepsilon}$ for any fixed $\varepsilon>0$ would be
very significant progress.

Our next question concerns forbidden minors.

{\bf Question 2.} Amongst all graphs $G$ not containing a $K_k$ minor
how large can $c(G)$ be?

As stated earlier Andreae \cite{Andreae} showed that $c(G) \leq {k-1 \choose 2}$ in this case.
Note that an upper bound that is less than quadratic in $k$ would be of great interest, because
if $c(G)=O(k^{2-\varepsilon})$ then it follows that
$c(n)=O(n^{1-\delta})$, where $\delta=\frac{\varepsilon}{4-\varepsilon}$. Indeed, if $G$ has a vertex with degree $\Omega(n^{\delta})$ then one cop can defend the neighbourhood of this vertex, so we proceed by induction. Else $G$ has $O(n^{1+\delta})$ edges, hence
the largest complete minor has at most $O(n^{1/2+\delta/2})$ vertices. Now $O(n^{(2-\varepsilon)(n^{1/2+\delta/2})}=O(n^{1-\delta})$ cops will suffice by our hypothesis.

Finally, for graphs on surfaces, Quilliot \cite{Q2} gave bounds of
$c(G) \leq 2k+3$ for a graph with orientable genus at most $k$. It
would be interesting to know what the true answer is.

{\bf Acknowledgement}. We would like to thank Jan Kratochvil for bringing this
problem to our attention, and for some interesting conversations.

\end{document}